# COMPARISON OF GEOMETRIC FIGURES

*Spyros Glenis M.Ed*

*University of Athens, Department of Mathematics,
e-mail spyros_glenis@sch.gr*

# Introduction

In Euclid, the geometric equality is based on the capability of superposition of the figures:

**Common notion 4**

*Things which coincide with one another are equal to one another.* ([12])

The geometric equality, with respect to Klein's view, is based on the group theory as well as on the set theory:

**Definition**

*Let a set $X \neq \varnothing$, G a subgroup of Aut(X) and the figures $S_1, S_2 \subseteq X$. We shall say that these figures are $G$–geometrically equal if and only if there is an $f \in G$ so that $f(S_1) = S_2$.*

The equality, indirectly defines the inequality of geometric figures. Euclid considers that a figure is smaller than another one if with an appropriate rigid motion the first coincides with part of the second. Although for any two figures $S_1, S_2$ it is easy to decide whether they are equal or not, however it is not that simple to decide if one of them is "smaller" than the other. Obviously a triangular region is never equal to a circular disk, but can we say that a triangular region is smaller than a circular disc if the radius of the disc is greater than or equal to the radius of the circumscribed circle of the triangle? In Euclid, the comparison involves only "similar" figures. On the contrary, Klein's view of equality, leads us to define a geometric inequality using the notion of "being subset" and enables us to compare even non-similar figures:

*We will say that $S_1$ is equal to or smaller than $S_2$ whenever there is a euclidean rigid motion $f$ so that $f(S_1) \subseteq S_2$. Then we will write $S_1 \leq S_2$.*

This "natural" definition of inequality provides a paradox as we will immediately illustrate using the following example given by the professor V. Nestorides:

Let us consider a closed half plane $A$ and let $B$ be the half plane $A$ with a line segment attached vertically to the edge of the half plane. Since $A \subseteq B$ we can say that $A \leq B$. Moreover, there is a translation of $B$, so that it is fully covered by $A$ and in this case we may write $B \leq A$. It seems logical to assume that $A$ and $B$ must be geometrically equal, in other words, that they can coincide if we apply a certain rigid motion. But this is impossible to happen, because every half plane remains half plane whenever we apply a rigid motion to it and obviously it can't coincide with a geometric figure that is not a half plane.

Since the geometric relation "$\leq$" is not antisymmetric it is necessary to restrict the comparison to certain classes of geometric figures. We already know that in the class of the line segments or in the class of the arcs of a circle, the relation "$\leq$" is a total order. Therefore the question is, if there are other classes of figures where the relation "$\leq$" is a total or a partial order.

We shall call *good classes* (of geometric figures) those that among the figures they contain we can't find a paradox like the one mentioned above. A good class, but not the only one, is that of the compact figures (sets). In fact, compact figures have the property not to generate paradox with any other geometric figure whether compact or not. Those figures will be called *good figures*. Besides the compacts, good figures are also the open-and-bounded sets. On the contrary, just bounded figures may not be good as we will prove later using a counterexample, given again by professor V. Nestorides.

The study, concerns not only the Euclidean Geometry, but it is also expanded into the Hyperbolic and the Elliptic Geometry and some parts may be formulated in a pure algebraic language so that they cover uniformly all three geometries. The conclusions we have reached, are fully compatible with our previous knowledge about the comparison of geometric figures. In the special case of the Euclidean Geometry we proved that there is a good class, containing all the fundamental geometric figures, where we can compare even non-similar ones. Therefore a comparison between a circular disc and a triangular region is meaningful in the new context.

# 1 Comparison of Figures in Euclidean Geometry

## 1.1 *Basic definitions*

We adopt Klein's view for the Euclidean Geometry. Our space is $\mathbb{R}^2$ endowed with the euclidean metric and the group acting on $\mathbb{R}^2$ is $ISO(\mathbb{R}^2)$, the group of euclidean isometries. The couple $(ISO(\mathbb{R}^2), \mathbb{R}^2)$ generates the euclidean geometric space where we will develop our study.

**Definition 1.1**

*Two figures $S_1$ and $S_2$ are geometrically equal when there is a euclidean isometry[1] $f: \mathbb{R}^2 \to \mathbb{R}^2$ so that $f(S_1) = S_2$. In that case we will write $S_1 \approx S_2$.*

Remarks

I. Figure is any subset of $\mathbb{R}^2$. From now on we will not distinguish the terms "subset of $\mathbb{R}^2$" and "figure".
II. We use the terms "rigid motion" and "isometry" synonymously.

**Definition 1.2**

*For any two figures $S_1$ and $S_2$ we shall say that $S_1$ is equal to or smaller than $S_2$ when there is a euclidean rigid motion $f$ so that $f(S_1) \subseteq S_2$. Then we will write $S_1 \leq S_2$.*

This "natural" definition does not satisfy in general the antisymmetric property, as we will prove later.

**Proposition 1.1**

*The relation $"\leq"$ is a pre-order of figures i.e. it is reflexive and transitive, and the reflex ion is meant in the sense of the geometric equality defined in 1.1*

Proof

---

[1] based on the euclidean metric $\rho$ of $\mathbb{R}^2$ where $\rho((x,y),(a,b)) = \sqrt{(x-a)^2 + (y-b)^2}$

Let A and B two geometrically equal figures. Then, by definition, there is a euclidean isometry $f$ so that $f(A) = B$. Then $f(A) \subseteq B$ also holds and we conclude that the relation $\leq$ is reflexive with respect to the geometric equality.

If $A \leq B$ and $B \leq C$ then there are isometries $f, g$ so that $f(A) \subseteq B$ and $g(B) \subseteq C$. Then for the isometry $g \circ f$ holds $g \circ f(A) \subseteq C$ i.e. $A \leq C$. Therefore the relation is transitive ▫

In the following examples we shall prove that "$\leq$" does not satisfy in general the antisymmetric property, with respect to the geometric equality of definition 1.1.

**Example 1.1**

Let the half lines $A = \{(x,0) \in \mathbb{R}^2 : x \geq 0\}$, $B = \{(x,0) \in \mathbb{R}^2 : x > 0\}$. Since $A$ is a closed subset of $\mathbb{R}^2$ while $B$ is not, **there is not** an isometry $f$ so that $f(A) = B$ [2] thus $A \not\approx B$.

For the isometry $f(x,y) = (x+1, y)$, $f(A) \subseteq B$ holds. But it is also obvious that $B \subseteq A$, so we have both $A \leq B$ and $B \leq A$ while $A \not\approx B$.

**Example 1.2**

Let the figure $A = \{(x,y) \in \mathbb{R}^2 : x \leq 0\} \cup \{(x,2) \in \mathbb{R}^2 : 0 \leq x \leq 1\}$ and the half plane $B = \{(x,y) \in \mathbb{R}^2 : x \leq 2\}$. Obviously $A \subseteq B$ therefore $A \leq B$

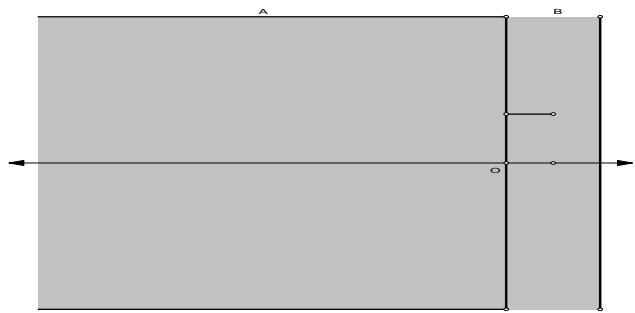

Since every isometry maps half planes into half planes there is not an isometry $f$ such that $f(A) = B$, so $A \not\approx B$.

---
[2] Every isometry maps closed sets into closed sets.

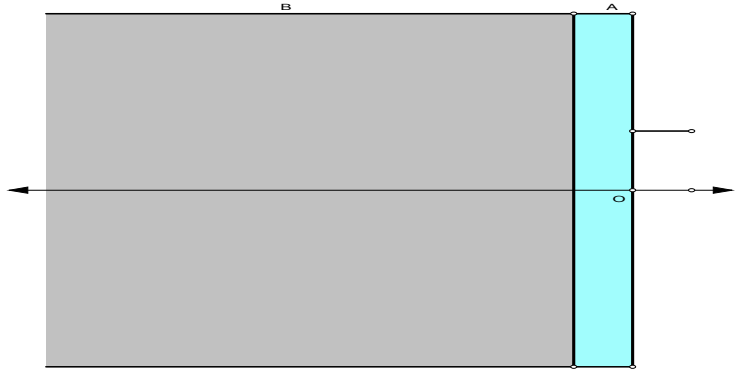

But for the isometry $f(x,y)=(x-3,y)$, $f(B) \subseteq A$ holds.

So we conclude that $A \leq B$ and $B \leq A$ while $A \not\cong B$.

**Example 1.3**

Let the figure $A = \{(x,y) \in \mathbb{R}^2 : x \geq 0, y \geq 0\}$ which is the right angle $\widehat{xOy}$ and the figure $B$ produced by $A$ when we subtract the inner part of the isosceles right triangle $S = \{(x,y) \in \mathbb{R}^2 : x \geq 0, y \geq 0, x+y < 1\}$.

Obviously $B \subseteq A$ so $B \leq A$ holds.

By translating $A$ parallel to the axis $x'x$ by two units we also have that $f(A) \subseteq B$, where $f(x,y)=(x+2,y)$ is an isometry. Thus $A \leq B$.

But $A$ and $B$ are not geometrically equal, because in case there is an isometry $g$ such that $g(A)=B$ and since isometries preserve the angles then $B$ should also be a right angle, which is absurd.

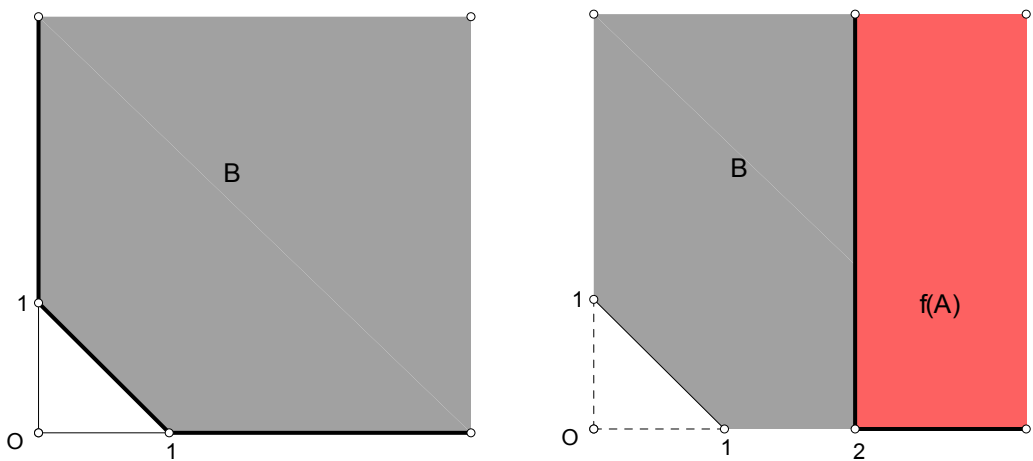

**Example 1.4**

Let an angle $\omega$ so that $\dfrac{\omega}{\pi}$ is irrational, for instance we can choose $\omega = \pi\sqrt{2}$ where $\sqrt{2}$ may be replaced by any irrational number. We set $a_k(\cos k\omega, \sin k\omega)$ a sequence of points lying on the circumference of the unit circle. By definition $a_k \neq a_m$ for $k \neq m$ otherwise we would have integers $n, \lambda$ so that $\dfrac{\omega}{\pi} = \dfrac{\lambda}{n}$.

The set $A = \{a_k : k = 0, 1, 2, 3, ...\}$ is a dense and equally distributed subset of the circumference. We also consider the set $B = A \setminus \{a_1\} \subseteq A$, therefore $B \leq A$. If there is an isometry $f$ of the plane such that $f(A) = B$ then we will arrive at a contradiction. For every $a \in A$ there is at least one $a' \in A$ so that the distance $d(a, a') = r$ where $r = 2\sin(\omega/2)$. Since $f(A) = B$ then and for every $b \in B$ there is at least one $b' \in B$ so that the distance $d(b, b') = r$. But this does not hold for $a_o \in B$ and we arrived at a contradiction.

Let $T$ be a rotation by $2\omega$ with center the origin of the axe. Then $T$ is an isometry and $T(A) = \{a_k : k = 2, 3, ...\} \subseteq B$. So $A \leq B$ and $B \leq A$ hold, while $A \not\approx B$.

We improve the definition of the pre-order "$\leq$" so that we arrive at an order relation:

**Definition 1.3**

*In the set of figures we define a relation $\lambda$ such that:*
$A \lambda B \Leftrightarrow \{A \approx B\} \text{ or } \{A \leq B \text{ and not } B \leq A\}$

**Proposition 1.2**

*$\lambda$ is an order relation.*

<u>Proof</u>

simple

**Definition 1.4**

*A class $\mathcal{E}$ of figures is said to be **good** when there are not any figures $A, B$ in the class $\mathcal{E}$ so that $A \leq B, B \leq A$ and $A \not\approx B$.*

<u>Remark</u>

The relation "$\lambda$" is defined so that it expels all the pathological cases of "$\leq$" where the antisymmetric property does not hold true. But the definition is not a "natural" one and gives no information to the question:

*Which figures form a good class?*

So the problem of the geometric order relation still remains unsolvable under the definition of $\lambda$. It is much wiser to stick to the definition of "$\leq$" and concentrate our study on those sets that satisfy the antisymmetric property.

**Definition 1.5**

*We will say that a figure $A$ is **good** when for every figure $B$, if $A \leq B$ and $B \leq A$ hold, then $A \approx B$ also holds.*

Obviously a class consisting only of good sets is a good class. The converse does not hold true. A trivial case is a class consisting of one figure (and all the geometric equals) that is not good. Since there is no other figure in the class to provide a counterexample then the class is good.

Another non-trivial case is the class of the open or closed angles. We proved in example 1.2 that an angle is not a good set but it is quite easy to verify that using open or closed angles only, we can not provide a counterexample.

## 1.2 *Quest for good classes of figures*

**Proposition 1.3**

*If $A$ is a good set and $A \approx B$ then $B$ is also a good set.*

<u>Proof</u>

Let $C \leq B$ and $B \leq C$. Since $A \approx B$ then $A \leq C$ and $C \leq A$ hold. As $A$ is a good set then there is an isometry $f$ of the plane such that $f(A) = C$. There is also an isometry $g$ of the plane such that $g(B) = A$. Then $g \circ f(B) = C$ i.e. $B \approx C$ □

**Proposition 1.4**

*If $A$ is a good set then its complement is also a good set.*

Proof

Let $B \le A^c$ and $A^c \le B$. Then there are isometries $f,g$, of the plane, so that $f(B) \subseteq A^c$ and $g(A^c) \subseteq B$. But then $f(B)^c \supseteq A$ and $g(A^c)^c \supseteq B^c$ hold.

Since $f, g$ are 1-1 and onto, $f(B)^c = f(B^c)$ and $g(A^c)^c = g(A)$ hold. Therefore $f(B^c) \supseteq A$ and $g(A) \supseteq B^c$ which is equivalent to $A \le B^c$ and $B^c \le A$. As $A$ is good there is an isometry $h$ of the plane such that $h(A) = B^c$. Then $h(A)^c = B$ and as $h$ is 1-1 and onto $h(A^c) = B$ i.e. $A^c \approx B$ □

**Theorem 1.1**

Let $\langle X, d \rangle$ be a compact metric space. If $f : X \to X$ is an isometry then $f(X) = X$.

Proof

Well known

**Proposition 1.5**

Every compact subset of $\mathbb{R}^2$ is a good set.

Proof

Let $A$ compact and $B$ an arbitrary set so that $A \le B$ and $B \le A$. Then there are isometries $f, g$ so that $f(A) \subseteq B$ and $g(B) \subseteq A$.

For the isometry $g \circ f : A \to A$ we have already seen that $g \circ f(A) = A$ because $A$ is compact (theorem 1.1).

But then $g(B) \subseteq A = g \circ f(A)$ which implies that $B \subseteq f(A)$. So $f(A) = B$ holds and $A \approx B$. Therefore $A$ is a good set □

Now we shall introduce a new definition that will be particularly useful.

**Definition 1.6**

*A figure $A \subseteq \mathbb{R}^2$ will be called **strongly good** if for every isometry $f : \mathbb{R}^2 \to \mathbb{R}^2$ that satisfies $f(A) \subseteq A$, then the equality $f(A) = A$ holds true.*

**Proposition 1.6**

*Every compact subset of $\mathbb{R}^2$ is strongly good.*

Proof

Direct from definition 1.6 and theorem 1.1 □

**Proposition 1.7**

*Every strongly good set is also a good set.*

Proof

Let $A$ a strongly good set and $B$ an arbitrary set such that $A \leq B$ and $B \leq A$. Then there are isometries $f, g$ of the plane so that $f(A) \subseteq B$ and $g(B) \subseteq A$.

For the isometry $g \circ f : A \to A$, $g \circ f(A) = A$ holds since $A$ is strongly good. Then $g(B) \subseteq A = g \circ f(A)$ therefore $B \subseteq f(A)$.

Finally we conclude that $f(A) = B$ so $A \approx B$ and $A$ is a good set □

Remark

The definition of the good figure is difficult to handle, as it depends on the "interaction" with all the other figures. On the contrary the definition of the strongly good figure is intrinsic, because, in simple words, strongly good is any figure that does not fit (without decomposition) into part of itself.

**Proposition 1.8**

*Every open and bounded subset of $\mathbb{R}^2$ is strongly good.*

Proof

Let $A$ open and bounded and $f$ an isometry of the plane such that $f(A) \subseteq A$. As $A$ is open then $A \cap \partial A = \emptyset$ and $\overline{A} = A \cup \partial A$. Thus $A = \overline{A} \setminus \partial A$

Also $diam(\overline{A}) = diam(A) < \infty$ so $\overline{A}$ is closed and bounded subset of $\mathbb{R}^2$, hence $\overline{A}$ is compact. Also $\partial A \subseteq \overline{A}$ and $diam(\partial A) \leq diam(\overline{A})$ so the boundary is closed and bounded hence compact subset of $\mathbb{R}^2$.

$f(A) \subseteq A \Rightarrow \overline{f(A)} \subseteq \overline{A} \Rightarrow f(\overline{A}) \subseteq \overline{A}$ and from T*heorem 3.1* $f(\overline{A}) = \overline{A}$ holds.

$f(\partial A) = f(\overline{A} \setminus A) = f(\overline{A}) \setminus f(A) = \overline{A} \setminus f(A) \supseteq \overline{A} \setminus A = \partial A$

Hence $f^{-1}(\partial A) \subseteq \partial A$ and as $f^{-1}$ is an isometry then from Theorem *3.1* again, we conclude that $f^{-1}(\partial A) = \partial A \Leftrightarrow \partial A = f(\partial A)$. Therefore $f(A) = f(\overline{A} \setminus \partial A) = f(\overline{A}) \setminus f(\partial A) = \overline{A} \setminus \partial A = A$, so $A$ is strongly good □

**Proposition 1.9**

*In $\mathbb{R}^2$, the union of a compact with an open bounded set is a strongly good set.*

Proof

Let $K$ compact and $A$ open and bounded subsets of $\mathbb{R}^2$, and the isometry $f$ of the plane such that $f(A \cup K) \subseteq A \cup K$.

The set $V = (A \cup K)^o \subseteq A \cup K$ is also open and bounded.

$W = K \setminus V = K \cap V^c$ is compact as an intersection of a compact with a closed set.
Obviously $V \cap W = \emptyset$

Since $A$ is open subset of $A \cup K$ then $A \subseteq V$ and

$A \cup K \subseteq V \cup K = V \cup (K \setminus V) = V \cup W \subseteq A \cup K$. Hence $V \cup W = A \cup K$

Also $V \cap W = \emptyset \Rightarrow f(V \cap W) = \emptyset \Rightarrow f(V) \cap f(W) = \emptyset$

Consequently $f(V \cup W) \subseteq V \cup W \Rightarrow f(V) \cup f(W) \subseteq V \cup W$

$f(V)$ is an open and bounded subset of $A \cup K$ therefore $f(V) \subseteq V$ since $V = (A \cup K)^o$. From *proposition 1.8* we conclude that $f(V) = V$.

As $f(V) \cup f(W) \subseteq V \cup W$ and $V \cap W = \emptyset, f(V) \cap f(W) = \emptyset$ then $f(W) \subseteq W$

But $W$ is compact and from *Theorem 3.1* $f(W) = W$. Therefore

$f(V) \cup f(W) = V \cup W \Rightarrow f(V \cup W) = A \cup K \Rightarrow f(A \cup K) = A \cup K$ □

**Proposition 1.10**

*In $\mathbb{R}^2$, the intersection of a compact with an open bounded set is a strongly good set.*

Proof

Let $K$ compact, $A$ open and bounded and $f$ an isometry of the plane such that $f(A \cap K) \subseteq A \cap K$

We set $X = \overline{A \cap K}$ which is a compact subset of $K$.

Then $f(A \cap K) \subseteq A \cap K \Rightarrow f(\overline{A \cap K}) \subseteq X \Rightarrow f(X) \subseteq X$. According to *Theorem 1.1* we conclude that $f(X) = X$.

Since $A\cap K$ is open in $K$, then it is also open in every closed subset of $K$. Therefore $A\cap K$ is open in $X$, so $X\setminus(A\cap K)$ is compact and it is obvious that $f^{-1}(X\setminus(A\cap K))\subseteq X\setminus(A\cap K)$. According to *Theorem 1.1* we conclude that $f^{-1}(X\setminus(A\cap K))=X\setminus(A\cap K)$.

But then $f(X\setminus(A\cap K))=X\setminus(A\cap K)$ and since $f(X)=X$ and $A\cap K\subseteq X$, we conclude that $f(A\cap K)=A\cap K$ i.e. $A\cap K$ is strongly good ☐

Remark

The proposition holds even if $A$ is not bounded.

**Proposition 1.11**

*The classes $\mathcal{X}=\{K\cup A$ where $K$ compact and $A$ open bounded subsets of the plane$\}$, $\mathcal{Y}=\{K\cap A$ where $K$ compact and $A$ open bounded subsets of the plane$\}$, $\mathcal{E}=\mathcal{X}\cup\mathcal{Y}$ and $\mathcal{F}=\mathcal{E}\cup\{S\subseteq\mathbb{R}^2:S^c\in E\}$ are all good classes.*

Proof

The classes above consist of strongly good sets. Then from definition 1.4 and the propositions 1.7, 1.8, 1.9, 1.10 the conclusion is obvious ☐

The class $\mathcal{F}$ includes almost all the fundamental figures of Euclidean Geometry: line segments, triangles, polygons, circles, arcs etc but not the open or closed angles which, as we have already mentioned, form a good class.

**Proposition 1.12**

If $\mathcal{W}$ is the class of open or closed angles then the class $\mathcal{F}\cup\mathcal{W}$ is good.

Proof

Since we already know that the classes $\mathcal{F},\mathcal{W}$ are good, then it is sufficient to examine whether a set from one class provides a counterexample to the other class.

Let $A$ a set of the class $\mathcal{F}$. Then $A\in\mathcal{E}$ or $A\in\mathcal{F}\setminus\mathcal{E}$.

If $A\in\mathcal{E}$ then it is bounded and it cannot provide a counter example with any angle of $\mathcal{W}$.

If $A\in\mathcal{F}\setminus\mathcal{E}$ and provides a counterexample with an angle then its complement $A^c\in\mathcal{E}$, will also provide a counterexample with the complement of the angle (which is also an angle). But this contradicts what we have already proved about the elements of $\mathcal{E}$ which give no counterexamples with the elements of $\mathcal{W}$.

So the class $\mathcal{F} \cup \mathcal{W}$ is a good one and includes the closed and open angles □

Remarks

1. The **closed sets** are neither good nor form a good class according to the example 1.1
2. The **connected sets** are neither good nor form a good class according to the example 1.1
3. The **bounded sets** are neither good nor form a good class according to the example 1.4
4. The **connected and bounded sets** are neither good nor form a good class. This can be easily proved if in the set $A$ of the example 1.4 we attach the inner points of the unit disc.
5. The **convex sets** are neither good nor form a good class according to the example 1.3
6. The **convex and bounded sets** are neither good nor form a good class. This can be easily proved if in the set $A$ of the example 1.4 we attach the inner points of the unit disc.

If two figures $A, B$ are good we have proved that their complements are also good. However the union or intersection of (strongly) good sets is not necessarily a good set as we will illustrate in the following counterexamples:

**Claim 1**

The figure $L = \{(x,0) \in \mathbb{R}^2 : x \geq 0\} \cup \{(0,1)\}$ is strongly good.

Proof

We denote $A = \{(x,0) : x \geq 0\} \subseteq \mathbb{R}^2$ and $b = \{(0,1)\}$ therefore $L = A \cup b$. From now on we will regard the point $(0,1)$ and the set $\{(0,1)\}$ as the same.

Let an isometry $f : \mathbb{R}^2 \to \mathbb{R}^2$ such that $f(L) \subseteq L \Leftrightarrow f(A \cup b) \subseteq A \cup b$

Since all isometries are 1-1 and $b \notin A$ then $f(b) \notin f(A)$.

The image of a half line through an isometry still remains a half line with point $f(a_o)$ as the origin, where $a_o = (0,0)$, thus $f(A) \subseteq A$

If $f(b) \in A$ then the points $f(b), f(a_o), f(a_1)$, where $a_1 = (5,0) \in A$ are collinear.

But then the points $b, a_o, a_1$ would also be collinear which is absurd!

Therefore $f(b) = b$.

Also $d(f(a_o), b) = d(f(a_o), f(b)) = d(a_o, b) = 1$

As $f(a_o) \in A$ and the only point $a \in A$ such that $d(a,b) = 1$ is $a_o$, then $f(a_o) = a_o$.

For every $a \in A$ we have that $d(a, a_o) = d(f(a), f(a_o)) = d(f(a), a_o)$

But there is only one point of the half line with this property, so $f(a) = a$.

We proved that $f = id_{\mathbb{R}^2}$ so for $L = A \cup b$ we have $f(L) = L$ □

**Claim 2**

The figure $M = \{(x,0) \in \mathbb{R}^2 : x \geq 0\} \cup \{(x,1) \in \mathbb{R}^2 : x < 0\}$ is strongly good.

Proof

We denote $B_1 = \{(x,0) : x \geq 0\}$, $B_2 = \{(x,1) : x < 0\}$. Then $M = B_1 \cup B_2 \subseteq \mathbb{R}^2$.

Let an isometry $f : \mathbb{R}^2 \to \mathbb{R}^2$ where $f(M) \subseteq M$.

$f(M) \subseteq M \Leftrightarrow f(B_1 \cup B_2) \subseteq B_1 \cup B_2 \Leftrightarrow f(B_1) \cup f(B_2) \subseteq B_1 \cup B_2$

Since $B_1 \cap B_2 = \emptyset \Leftrightarrow f(B_1 \cap B_2) = \emptyset \Leftrightarrow f(B_1) \cap f(B_2) = \emptyset$

The images of $B_1, B_2$ are half lines therefore $f(B_1) \subseteq B_1$ and subsequently $f(B_2) \subseteq B_2$ or $f(B_1) \subseteq B_2$ and subsequently $f(B_1) \subseteq B_2$.

If $f(B_1) \subseteq B_1$ and $f(B_2) \subseteq B_2$ then $f(0,0) = (x_1, 0)$ for some $x_1 \geq 0$ and

$f((0,1)) = (x_2, 1)$ for some $x_2 \leq 0$

$$1 = d((0,0),(0,1)) = d(f(x_1,0), f(x_2,1)) = d((x_1,0),(x_2,1)) = \sqrt{(x_1-x)^2+1}$$

Then $|x_1 - x_2| = 0$ and we conclude that $x_1 = x_2 = 0$.

Therefore $f((0,0)) = (0,0)$ and $f((0,1)) = (0,1)$

If $b_o = (0,0)$ and $b_1 = (1,0)$ then

$$1 = d(b_o, b_1) = d(f(b_o), f(b_1)) = d((0,0), f(b_1))$$

As $f(b_1) \in f(B_1) \subseteq B_1$ then $f(b_1) \equiv b_1$

Since $f$ preserves the three $(0,0)$, $(0,1)$ and $(1,0)$ which are non-collinear then

$$f = id_{\mathbb{R}^2}$$

If $f(B_1) \subseteq B_2$ and $f(B_2) \subseteq B_1$ then $f(0,0) = (x_1, 1)$, $x_1 < 0$ and $f(0,1) = (x_2, 0)$, $x_2 \geq 0$.

But $1 = d((0,0),(0,1)) = d((x_1,1),(x_2,0)) = \sqrt{(x_1-x_2)^2+1}$ then $x_1 = x_2 = 0$ which is absurd as $x_1 < 0$.

So the only case is $f = id_{\mathbb{R}^2}$ and then $f(M) = M$ □

**Example 1.5**

Using the previous notations for the sets $M$ and $L$, then the set $L \cap M = \{(x,0) \in \mathbb{R}^2 : x \geq 0\}$ is not good according to the example 1.1, although $L, M$ are strongly good sets.

**Example 1.6**

We also use here the previously defined sets $L$ and $M$.

If $g$ is a reflection with respect to the $y'y$ axis then $g(M)$ is strongly good but the set $G = L \cup g(M) = \{(x,0) \in \mathbb{R}^2 : x \in \mathbb{R}\} \cup \{(x,1) \in \mathbb{R}^2 : x \geq 0\}$ is not good because for the set $V = \{(x,0) \in \mathbb{R}^2 : x \in \mathbb{R}\} \cup \{(x,1) \in \mathbb{R}^2 : x > 0\}$ we have that $V \subseteq G$ and $h(G) \subseteq V$, where $h(x, y) = (x+1, y)$. Therefore $V \leq G$ and $G \leq V$. But $G$ is a closed set and $V$ is not closed so there is no isometry $f$ such that $V = f(G)$, hence $V \napprox G$.

# 2  Good figures and good classes in non‑Euclidean Geometries

## 2.1  *Hyperbolic Geometry*

We will use the *Poincare's* disc as model. The non-euclidean plane is the interior of the unit disc $D$ of the complex plane i.e.

$$D = \{z \in \mathbb{C} : |z| < 1\}$$

The distance between the points $A(z_1), B(z_2)$ of $D$ is given by the formula

$$d_p(A, B) = 2\operatorname{arctanh} \frac{|z_2 - z_1|}{|1 - \overline{z}_1 z_2|} \quad (1)$$

and we will call it the *hyperbolic metric*.

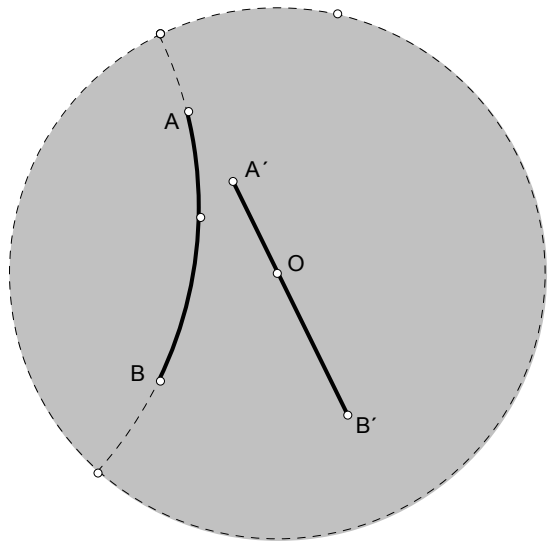

The isometries of the model are all functions $f : D \to D$ of the form:

$$f(z) = \frac{az + b}{\overline{b}z + \overline{a}} \text{ and } a\overline{a} - b\overline{b} = 1 \text{ (without change of orientation) or}$$

$$f(z) = \frac{a\overline{z} - b}{\overline{b}\overline{z} - \overline{a}} \text{ και } a\overline{a} - b\overline{b} = 1 \text{ (with change of orientation)}$$

**Definition 2.1**

*Two figures $S_1$ and $S_2$ of the hyperbolic plane D, are geometrically equal whenever there is an isometry[3] $f: D \to D$ so that $f(S_1) = S_2$. Then we will write $S_1 \approx S_2$*

**Definition 2.2**

*For two figures $S_1$ and $S_2$ of the hyperbolic plane D, we will say that $S_1$ is equal to or smaller than $S_2$ whenever there is an isometry $f$ such that $f(S_1) \subseteq S_2$. Then we will write $S_1 \leq S_2$.*

The definitions of a **good set**, **strongly good set** and **good class** in the Hyperbolic Geometry are similar to those of the Euclidean Geometry.

The relation "$\leq$" is reflexive, transitive but not antisymmetric as we will prove using the following examples.

**Example 2.1**

The example 1.4 of the irrational rotation $\omega$ in the euclidean plane, is also applied in the hyperbolic plane provided that the center of rotation is the origin and the points $a_k$ are $a_k = (b\cos k\omega, b\sin k\omega)$ with $0 < b < 1$. It is noted here that the angles of the hyperbolic plane $D$ are the same with the corresponding euclidean ones.

It can be proved, exactly as in the euclidean plane, that the sets $A = \{a_k : k = 0,1,2,3,...\}$ and $B = A \setminus \{a_1\}$ are not geometrically equal.

It is also true that $B \leq A$ because $B \subseteq A$.

The rotation $T(z) = e^{2i\omega}z$, is an isometry of the hyperbolic plane $D$ because it is the case of $f(z) = \dfrac{az+b}{\bar{b}z+\bar{a}}$ with $b = 0$, $a = e^{i\omega}$ and $a\bar{a} - b\bar{b} = 1$.

Performing the operations we reach that $T(a_k) = a_{k+2}$ so

$T(A) = \{a_k : k = 2,3,...\} \subseteq B$ i.e. $A \leq B$.

Since $A \leq B$, $B \leq A$ while $A \not\approx B$ then $A$ is not a good set in Hyperbolic Geometry.

---

[3] with respect to the hyperbolic metric $d_p$

**Example 2.2**

Another model for the Hyperbolic Geometry is given by the half plane
$$H = \{z \in \mathbb{C} : \operatorname{Im}(z) > 0\}$$

In $H$ the metric is given by
$$d_H(z_1, z_2) = \operatorname{arctanh} \frac{|z_2 - z_1|}{|\bar{z}_2 - z_1|}$$

We point out that the models of the disc and of the half plane are isometrically isomorphic and the corresponding mapping is:
$$h : H \to D \;\; \mu\varepsilon \;\; h(z) = \frac{z - i}{z + i}$$

We consider $A = \{z \in H : \operatorname{Re}(z) \geq 0\} \cup \left\{z \in H : |z| = 1 \;\mu\varepsilon\; \operatorname{Re}(z) \leq 0 \;\text{και}\; \operatorname{Im}(z) \geq \frac{1}{2}\right\}$

and $B = \{z \in H : \operatorname{Re}(z) \geq 1\}$. Obviously $B \subseteq A$ (see figure 4.2).

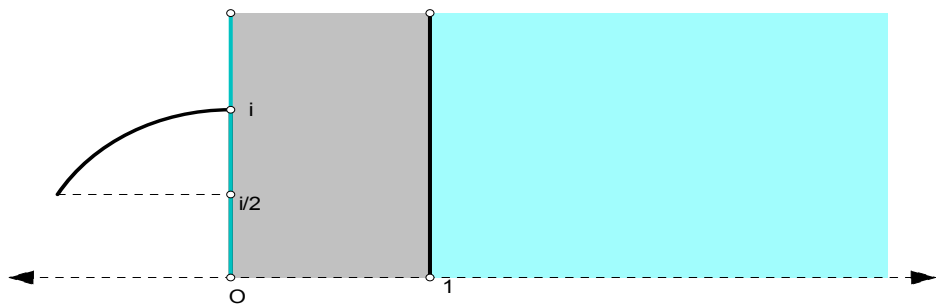

It is easy to verify that the mapping $g_H(z) = z - 3$ is an isometry in $H$.

Then $A$ is a subset of $g_H(B) = \{z \in H : \operatorname{Re}(z) \geq -2\}$.

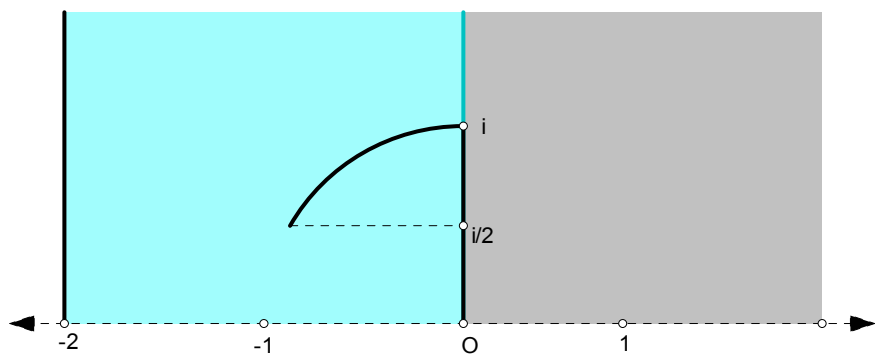

Consequently we have that $A \leq B$ and $B \leq A$. But there is no isometry $f$ such that $f(A) = B$ since in B every point has infinitely many points within a given hyperbolic distance $\rho > 0$; however in $A$ the point $-1 + \frac{i}{2}$ for small $\rho > 0$ has only one point in $A$ at hyperbolic distance $\rho$.

The Moebius transformation $w = \frac{z-i}{z+i}$ with $\text{Re}(z) > 0$ maps conformally the hyperbolic plane $H$ to the hyperbolic plane $D$ and transforms the sets $A, B$ and the isometry $g_H$ of $H$ into the sets $A', B'$ and the isometry $g_D$ of $D$ so that:

$$A' \leq B', \ B' \leq A' \text{ and } A' \not\approx B'.$$

Now we will search for good sets and good classes of the hyperbolic plane $D$

**Lemma 2.1**

*If $z, w \in D$ then $|1 - \bar{z}w| > 1 - |w|$*

Proof

Well known exercise.

**Proposition 2.1**

*The euclidean and the hyperbolic metric on the unit disc $D$ are equivalent.*

Proof

Well known

**Proposition 2.2**

*Let $A \subseteq D$ then $A$ is bounded with respect to the euclidean metric. However $D$ itself is bounded in euclidean metric but unbounded in hyperbolic metric.*

Proof

Obvious.

**Corollary 2.1**

*If $A$ is a closed and bounded subset of $D$ with respect to the hyperbolic metric then it is compact.*

Proof

Since the metrics are equivalent and $A$ is closed with respect to the hyperbolic metric then it is also closed (and bounded) with respect to the euclidean metric, hence compact. As the metrics are equivalent then it is also compact with respect to the hyperbolic metric ∎

**Proposition 2.3**

*In $D$, every compact and every open and bounded (in hyperbolic distance) set is a strongly good set. In addition the union, the intersection and the complements ot them are also strongly good sets of the Hyperbolic Geometry.*

Proof

The proofs are similar to those of proposition 1.6

**Proposition 2.4**

*If we are studying only the figures in the interior of the unit disc, then the class $\mathcal{F} \cup \mathcal{W}$ of proposition 1.12 is also a good class of the hyperbolic plane.*

Proof

It is obvious, by the definitions and the proposition 2.3.

## 2.2 Elliptic Geometry

As model of the Elliptic Geometry we will use the surface $S^+$ of the upper hemisphere of the Riemannian sphere i.e.:

$$S^+ = \{(x, y, z) \in \mathbb{R}^3 : x^2 + y^2 + z^2 = 1,\ z \geq 0\}$$

with the assumption that two antipodal points of the equator coincide in other words the points $a(x, y, 0) \in S^+$ and $a(-x, -y, 0) \in S^+$ make up of the same point $\bar{a}$.

The lines of the model are the intersections of $S^+$ with planes that pass through the center of the sphere. Then we can define the distance between two points of the elliptic plane as:

If $a, b \in S^+$ then

$$d_s(a, b) = \arccos|\vec{a} \cdot \vec{b}|$$

where $\vec{a} \cdot \vec{b}$ is the ordinary dot product of $\mathbb{R}^3$. We will call $d_s$ the *elliptic metric*.

**Definition 2.3**

*Two figures $S_1$ and $S_2$ of the elliptic plane $S^+$, are geometrically equal whenever there is an isometry[4] $f : S^+ \to S^+$ so that $f(S_1) = S_2$. Then we will write $S_1 \approx S_2$*

**Definition 2.4**

*For two figures $S_1$ and $S_2$ of the elliptic plane $S^+$, we will say that $S_1$ is equal to or smaller than $S_2$ whenever there is an isometry $f$ such that $f(S_1) \subseteq S_2$. Then we write $S_1 \leq S_2$.*

The definitions of a **good set, strongly good set** and **good class** of the Elliptic Geometry are similar to those of Euclidean and Hyperbolic Geometry.

The relation "$\leq$" is reflexive and transitive but not antisymmetric, as we will prove next.

Before proving that "$\leq$" is not antisymmetric we must note that:

α) The lines of the elliptic plane have finite length $\pi$

β) The elliptic plane is not oriented, in fact it is a one-sided closed surface (such a surface is the Moebius bundle). As such it is meaningless to talk about half planes or half lines and anything else defined by them. Therefore, the examples 1.1, 1.2, and 1.3 cannot be realized in the Elliptic Geometry.

**Claim**

*The mapping $T : S^+ \to S^+$ where $T(x,y,z) = \begin{pmatrix} \cos 2\omega & -\sin 2\omega & 0 \\ \sin 2\omega & \cos 2\omega & 0 \\ 0 & 0 & 1 \end{pmatrix} \begin{pmatrix} x \\ y \\ z \end{pmatrix}$ is an isometry of the elliptic plane $S^+$.*

Proof

Let $a(x_1, y_1, z_1)$ and $b(x_2, y_2, z_2)$ then

$T(\vec{a}) = (x_1 \cos\omega - y_1 \sin\omega, x_1 \sin\omega + y_1 \cos\omega, z_1)$,

$T(\vec{b}) = (x_2 \cos\omega - y_2 \sin\omega, x_2 \sin\omega + y_2 \cos\omega, z_2)$

---

[4] with respect to the elliptic metric $d_s$

and the dot product $T(\vec{a}) \cdot T(\vec{b}) = x_1 x_2 + y_1 y_2 + z_1 z_2 = \vec{a} \cdot \vec{b}$

So $d_s(T(a), T(b)) = \arccos|T(\vec{a}) \cdot T(\vec{b})| = \arccos|\vec{a} \cdot \vec{b}| = d_s(a,b)$

This isometry corresponds to a rotation of the elliptic plane by angle $2\omega$ with center the north pole $N(0,0,1)$ □

**Example 2.3**

We first consider some $\omega \in \mathbb{R}$ such that $\dfrac{\omega}{\pi}$ is irrational.

Let the points $a_k = \left(\dfrac{1}{2}\cos k\omega, \dfrac{1}{2}\sin k\omega, \dfrac{\sqrt{3}}{2}\right)$ of the surface $S^+$.

It is easy to verify that $\vec{a}_k \vec{a}_{k+1} = \dfrac{1}{4}\cos\omega + \dfrac{3}{4}$ is independent of $k$.

Then $d_s(a_k, a_{k+1}) = r$ is constant for every $k = 0, 1, 2, \ldots$

We consider the figures $A = \{a_k : k = 0, 1, 2, 3, \ldots\}$ and $B = A \setminus \{a_1\}$.

Since $B \subseteq A$ then $B \leq A$

Using the above isometry $T$ it is easy to see that $T(a_k) = a_{k+2}$ so $T(A) = \{a_k : k = 2, 3, \ldots\} \subseteq B$ i.e. $A \leq B$.

If there was an isometry of the elliptic plane $f : S^+ \to S^+$ so that $f(A) = B$ then we result in absurd as in the example 1.4.

So the relation "$\leq$" is not an order in the Elliptic Geometry.

**Proposition 2.5**

$S^+$ is compact (with respect to the elliptic metric)

<u>Proof</u>

Let the sequence $\{u_n\} \subseteq S^+$.

$\{u_n\}$ is bounded in $\mathbb{R}^3$ with respect to the euclidean metric, so it has a converging subsequence $\{u_{k_n}\}$ where $u_{k_n} \to u \Leftrightarrow (x_{k_n}, y_{k_n}, z_{k_n}) \to (x, y, z)$, and $u \in S^+$ since it is a closed subset of $\mathbb{R}^3$.

Then $d_s(u_{k_n}, u) = \arccos|\vec{u}_{k_n} \cdot \vec{u}| \to \arccos|u \cdot u| = \arccos 1 = 0$

We proved that every (bounded) sequence of $S^+$ has a converging subsequence therefore $S^+$ is compact ☐

**Corollary 2.2**

*Every closed subset of $S^+$ is compact.*

**Proposition 2.6**

*In $S^+$, every closed and every open set is a strongly good set of the Elliptic Geometry.*
Proof
Every closed set is compact in $S^+$, so it is also strongly good. Then its complement it is also strongly good, therefore and every open set is strongly good ☐

**Proposition 2.7**

*If A is open set and K is closed set of the elliptic plane then the set $A \cup K$ is strongly good.*
Proof
The same as in proposition 1.9

**Corollary 2.3**

*If A is open set and K is closed set of the elliptic plane then the set $A \cap K$ is strongly good.*
Proof
The set $A^c \cup K^c$ is strongly good according to the proposition 2.7
Then its complement $\left(A^c \cup K^c\right)^c = A \cap K$ is also strongly good ☐

## 3  Open issues

A fundamental question is whether the definitions of the good set and the strongly good set are equivalent or there is a counterexample of a good set that is not strongly good. In the appendix it is proved that in $\mathbb{R}$ all good sets are also strongly good. So it is our belief that the definitions are also equivalent on the plane.

Another question is whether the algebra produced by $X = \{A \subseteq \mathbb{R}^2 : A$ compact or open and bounded$\}$ consists only of strongly good sets.

Finally, as in $\Omega = \{$set of good classes$\}$ the assumptions of Zorn's[5] lemma are satisfied it would be quite interesting to find maximal classes within the set $\Omega$.

---

[5] Every totally ordered subset of $\Omega$ is defined to be a set of good classes $Y = \{F_i, i \in I\}$ so that for every $i, j \in I$ it is true that $F_i \subseteq F_j$ or $F_j \subseteq F_i$. Obviously the **good** class $\bigcup_{i \in I} F_i$ is an upper bound of $Y$.

# Appendix

In the appendix we prove that good sets coincide with strongly good sets in $\mathbb{R}^1$. We do not know what happens in $\mathbb{R}^2$.

**Lemma**

Let a set $A \subseteq \mathbb{R}$, $a \in A$ and an isometry $f : \mathbb{R} \to \mathbb{R}$ such that $f(A) = A \setminus \{a\}$ i.e. $A \approx A \setminus \{a\}$. Then $A$ is not a good set.

**Proof**

If $f : \mathbb{R} \to \mathbb{R}$ is an isometry then $f(x) = x$ or $f(x) = -x$ or $f(x) = -x + c$, $c \neq 0$ or $f(x) = x + c$, $c \neq 0$.

- If $f(x) = x$ then $f(A) = A \setminus \{a\} \Leftrightarrow A = A \setminus \{a\}$ i.e. $a \notin A$, absurd.

- Both $f(x) = -x$ and $f(x) = -x + c$ have the property $f^2 = id_{\mathbb{R}}$.

  From the assumption we have that $f(A) = A \setminus \{a\}$ therefore

  $f^2(A) = A \setminus \{a, f(a)\} \Leftrightarrow A = A \setminus \{a, f(a)\}$ i.e. $a \notin A$, absurd.

- If $f(x) = x + c$, $c \neq 0$.

  Then $f(a - c) = a$, so $a - c \notin A$

  Also $f(a) = a + c \in f(A) = A \setminus \{a\}$.

  Let the set $B = A \setminus \{a + c\}$. Then $a + c \notin B$ and $a - c \notin B$ since $a - c \notin A$ and $B \subseteq A$.

  Obviously $B \subseteq A$ and $f^2(A) = A \setminus \{a, a + c\} \subseteq B$

  Let us assume that there is an isometry $g : \mathbb{R} \to \mathbb{R}$ so that $g(A) = B$

  For every $x \in A$ we have that $f(x) = x + c \in A \setminus \{a\}$ therefore there is at least one $x' = x'(x) \in A$ so that $d(x, x') = c$ (for instance $x' = x + c$).

  But then we will also have that $d(g(x), g(x')) = c$.

  Since $g(A) = B$ there is some $x_o \in A$ such that $g(x_o) = a \in B$.

  Then there is also some $x'' = x''(x_o) \in A$ so that $d(x_o, x'') = d(x_o, x'(x_o)) = c$

It follows that $c = d(x_o, x'') = d(g(x_o), g(x'')) = d(a, g(x''))$.

We may say that there is some $b \in B = g(A)$, where $b = g(x'')$ so that

$d(a,b) = c$.

But then $b = a - c$ or $b = a + c$, which in either case are not points of $B = g(A)$.

This is a contradiction and we conclude that $A$ is not a good set.

**Proposition**

Let $A \subseteq \mathbb{R}$ be a good set, then $A$ is strongly good.

**Proof**

We assume that $A$ is not strongly good.

Then there will be an isometry $T: \mathbb{R} \to \mathbb{R}$ so that $T(A) \subsetneq A$.

Therefore there is $a \in A \setminus T(A)$ and $T(A) \subseteq A \setminus \{a\} \subseteq A$ i.e. $T(A) \leq A$

Since $A \subseteq T^{-1}(A) \Leftrightarrow A \subseteq T^{-1}(T^{-1}(T(A)))$ then $A \leq T(A)$ and as $A$ is a good set then $A \approx T(A)$.

We can also prove that $A \approx A \setminus \{a\}$ because:

$A \setminus \{a\} \leq A$ and $A \approx T(A) \leq A \setminus \{a\}$ since $A$ is a good set we conclude that $A \setminus \{a\} \approx A$.

But from the previous lemma such a set $A$ is never a good set, absurd! □